\documentclass[12pt]{article}
\usepackage{amsmath,amsthm,amsfonts,fancyhdr,placeins}
\usepackage{graphicx,enumerate,relsize}
\newcommand{\M}{m}

\newcommand{\g}{\lambda}

\newtheorem{thm0*}{Proposition}
\newtheorem{corollary}{Corollary}
\newtheorem*{thm*}{ Theorem, Nikishin }
\newtheorem*{thm2*}{ Theorem, Kostyuchenko-Mirzoev, 1998}
\newtheorem*{thm3*} {Theorem, Berg-Dur\'an, 2006}
\newtheorem*{thm4*}{Theorem 3.3 of \cite{osc}}
\newtheorem*{thm5*}{Lemma 3.2 of \cite{osc}}
\newtheorem{theorem}{Theorem}
\newtheorem*{corollary*}{Corollary}%[theorem]
\begin{document}
\title {On the finiteness of logarithmic Hamiltonians for Volterra-type lattices in terms of the spectral measures of Jacobi operators}
\author{ Andrey Osipov  \\
Scientific Research Institute for System Analysis \\
of the National Research Centre  ``Kurchatov Institute"\\
(NRC ``Kurchatov Institute" - SRISA),\\
Nakhimovskii pr. 36-1, Moscow, 117218, Russia \\  email:
{osipa68@yahoo.com} }

\date{}
\fancyhead{} \chead{}
 \maketitle

\begin{abstract}
We establish a correspondence between the semi-infinite and
infinite Volterra lattices having a finite logarithmic Hamiltonian
and certain classes of even probability measures. In doing so, we
apply the inverse spectral theory of Jacobi operators and the
theory of orthogonal polynomials. A similar correspondence is
established for semi-infinite modified Volterra lattices.
\end{abstract}


%%%%%%%%%%%%%%%%%%%%%%%%%%%%%%%%%%%%%%%%%%%%%%%%%%%%%%%%%%%%%%%%%%%%%%
%%%%%%%%%%%%%%%%%%%%%%%%%%%%%%%%%%%%%%%%%%%%%%%%%%%%%%%%%%%%%%%%%%%%%%
\section{Introduction}

Since the famous work of Gardner, Greene, Kruskal and Miura, the
integration of nonlinear equations by using the existing inverse
problem methods for differential and difference operators as well
as the development of these methods aiming to utilize them for
this integration task, is among the major topics of modern
mathematical physics. This trend found its reflection in the work
of Nikishin \cite{nik}, where the inverse scattering problem for
second-order difference operators (discrete Sturm-Liouville or
Jacobi operators) was studied and applied to the integration of
two well-known nonlinear dynamical systems: the Toda lattice and
the Volterra lattice (also known as Langmuir chain, Kac - van
Moerbeke system and discrete Korteweg - de Vries equation). The
lattices were considered both in the semi-infinite case (i. e.
then their discrete ``spatial" variable takes nonnegative integer
values, see next section) and the infinite case (for the discrete
variable running through all the integers). Instead of solving the
inverse scattering problem by using a discrete version of
Gel'fand-Levitan-Marchenko equation \cite{to,gus1,gus2,ser},
Nikishin proposed in \cite{nik} an approach based on the theory of
orthogonal polynomials on the unit circle and real line (the
results on their asymptotic properties \cite{seg,ger}), the links
between the latter and the spectral theory of Jacobi operators and
the expansion into continued $J-$ fractions of their Weyl
functions \cite{ah,nisor,bk, yurm, osj} (more details will be
given in the next section). He also studied the connection between
the inverse scattering problem and the inverse spectral problem
for Jacobi operators (well familiar in connection with the
Hamburger moment problem \cite{ah,ber,nisor,seg}) , which amounts
to reconstruction of the operator from its spectral measure. In
particular, this study involved determination of the spectral
function from the scattering data. This approach is applicable to
a wider class of scattering data (in \cite{nik} they were called
generalized scattering data) than in the aforementioned works
\cite{gus1,gus2,ser}. Using the Lax representation for both
lattices, the evolution in time of the corresponding generalized
scattering data was obtained, and the integration problem for the
lattices was reduced to a certain problem of complex analytic
factorization. Note that approximately at the same time when
\cite{nik} was published, Berezansky in \cite{bert1,bert2}
proposed an integration method for semi-infinite Toda lattice
based on the inverse spectral problem for Jacobi operators
\cite{ber} (he solved the initial boundary value problem for this
system in the class of bounded solutions). These works have much
in common in regards to ``integration" part. Since the scattering
data can be defined only for a certain class of bounded Jacobi
operators, the latter integration method is applicable to a wider
family of semi-infinite Toda lattices than the one of \cite{nik}.
At the same time, \cite{bert1,bert2} contain no study of the
spectral function (e. g., how to get one from the initial data)
assuming that it is given a priori. Overall, it can be said that
in \cite{nik}, a function-theoretical approach to the integration
of nonlinear lattices have been proposed, and its elements (e. g.
the application of functional continued fractions) were used in a
number of subsequent works \cite{sor,soriz,viz,apt}.

Here we will use the results of \cite{nik} in the study of
Hamiltonian properties (i.e. Hamiltonians, Poisson brackets,
master symmetries, etc.) of Volterra lattice and its modified
version. Such properties for Volterra-type systems (including
Bogoyavlensky lattice or hungry Volterra model
\cite{hik1,pan})have been subject of intense studies in recent
decades, see \cite{ft1,dam,sur,hik1,osajm}. The fact that Volterra
lattice is integrable by means of the inverse spectral problem for
Jacobi operators in particular implies that the Hamiltonians of
the lattice can be expressed via the spectral function of Jacobi
operator which appear in the Lax representation for the lattice.
As we will see below, the findings contained in \cite{nik} and
their further development, see \cite[Chaper 13]{simop} and
references thereafter, can be utilized to describe in terms of the
spectral measure the physically meaningful situation when the
logarithmic (or zeroth \cite{osrjnd}) Hamiltonian of Volterra
lattice is finite. It will be shown that for semi-infinite
Volterra lattices with the finite logarithmic Hamiltomian a
correspondence with a class of even spectral measures satisfying
Szeg\H{o} condition (or having finite entropy \cite{besden}) can
be established. Such correspondence in particular implies that
these spectral measures can be characterized by using the
logarithmic Hamiltonians. Note that a similar study of the
correspondence between canonical Hamiltonian systems on the
half-line and their spectral measures have been done by Bessonov
and Denisov in \cite{besden}, where a criterion of convergence of
logarithmic integral for a spectral measure of the Krein string
have been obtained. In the infinite case, the inverse spectral
problem method for Jacobi operators can be also applicable, in
particular, by using the ``duplication" procedure, see
\cite{bg,zh} and the next section. Here we will be able as well to
find a correspondence between the finite logarithmic Hamiltonians
of infinite Volterra lattices and certain couples of spectral
functions belonging to the same class as in the semi-infinite
case. At the end of the paper, we apply our results concerning
semi-infinite Voltera lattices to the study of modified Volterra
lattices \cite{sur} with a finite logarithmic Hamiltonian.

%%%%%%%%%%%%%%%%%%%%%%%%%%%%%%%%%%%%%%%%%%%%%%%%%%%%%%%%%%%%%%%%
%%%%%%%%%%%%%%%%%%%%%%%%%%%%%%%%%%%%%%%%%%%%%%%%%%%%
\section{Inverse spectral problem for Volterra lattice and its
Hamiltonian structure} We start with the classical Volterra (or
Kac - van Moerbeke) lattice
\begin{eqnarray}
\label{vl} \overset \cdot a_n = a_n(a_{n+1}-a_{n-1}),\qquad \quad\\
\quad \; a_n=a_n(t)>0,\quad t \in [0,T),\quad 0<T \le \infty;
\nonumber
\end{eqnarray}
in the semi-infinite case:
\begin{equation}
\label{seminf}
 i \in \mathbb{Z}_+,\quad \quad a_{-1}=0,
\end{equation}
with the initial data $a(0):=(a_n(0))_{n=0}^{\infty}.$ It admits
the Lax representation $\overset \cdot L = [L,A]:=LA-AL$ with the
following infinite matrices $L$ and $A$:
\begin{equation}
L=L(t)=(L_{i,j}(t))_{i,j=0}^{\infty}=\begin{pmatrix} 0& \sqrt{a_0}& 0&0&0&\dots\\
           \sqrt{a_0}&0&\sqrt{a_1}&0&0&\dots\\
           0&\sqrt{a_1}&0&\sqrt{a_2}&0&\dots\\
            0&0&\ddots&\ddots&\ddots&\dots
            \end{pmatrix};
\label{lv}
 \end{equation}
 \begin{eqnarray}
\label{av}
 &&A=A(t)=(A_{i,j}(t))_{i,j=0}^{\infty}=\qquad \qquad \qquad \qquad
 \qquad\nonumber
\\&& \nonumber
 \\&&=\frac{1}{2}\begin{pmatrix} 0& 0& -\sqrt{a_0a_1}&0&0&0&\dots\\
           0&0&0&-\sqrt{a_1a_2}&0&0&\dots\\
           \sqrt{a_0a_1}&0&0&0&-\sqrt{a_2a_3}&0&\dots\\
            0&\sqrt{a_1a_2}&0&0&0&-\sqrt{a_3a_4}&\dots\\
            0&0&\ddots&\ddots&\ddots&\ddots&\ddots
            \end{pmatrix}.\\\nonumber
\end{eqnarray}
Here the matrix $L$ is a Jacoby matrix with zero main diagonal
(sparse Jacobi matrix), $A$ is skew-symmetric. Denote by $\cal L$
the minimal closed symmetric operator generated by $L$ and acting
in the Hilbert space of complex quadratic summable sequences
$l^2[0,\infty)$ (we denote by $\{e_n\}_{n=0}^\infty$ its standard
orthonormal basis), see e. g. \cite{ah,ber}. Let $F_\g $ be a
spectral function of $\cal L$ (if $\cal L$ is non-selfadjoint,
then $F_\g$ is not unique), see \cite[Chapter IV]{ah}-\cite{agl}.
Then, set $\rho(\g)=(F_\g e_0,e_0)$ and define the moments: $s_k =
\int_{-\infty}^{\infty} \g^k d \rho(\g), \; k \in \mathbb{Z_+}$.
Since $L$ is a sparse matrix, it turns out that $d\rho(\g)$ is an
even measure and therefore all its moments $s_k$ for odd $k$ are
zero (more details about sparse Jacobi and band matrices can be
found in \cite{osrjmp}). The elements of $L$ can be reconstructed
from $d\rho(\g)$ as follows:
\begin{equation}
\label{solisp} a_n=\frac{\Delta_{n+1}\Delta_{n-1}}{\Delta^2_n},
\quad n\in \mathbb{Z_+},
\end{equation}
where $\Delta_{n}$ are determinants of the Hankel matrices $ H_n
=(s_{i+j})_{i,j=0}^n,\,$ $\Delta_{-1}=1,$ and by the Hamburger
theorem, $\Delta_n>0$ for all $n$, see \cite[Chapters 1,4]{ah}.

From now on, as in \cite{nik}, we will assume that the operator
$\cal L$ is bounded, i. e. $\sup_n a_n < \infty.$ Then it is
selfadjoint, see \cite[Chapter VII]{ber} or \cite[Chapter I]{ah}.
Let $\{E_\g\}$ be the corresponding spectral family/resolution of
the identity. Then
\begin{equation}
\label{dro} d\rho(\g)=d(E_\g e_0,e_0)
\end{equation}
is the spectral measure of $\cal L,$ so $ d\rho(\g) $ is a
probability measure on $\mathbb{R}$. Denote by $D({\cal L})$ the
domain of $\cal L$. Since for all $k,$ ${\cal L}^k e_0 \in D({\cal
L})$ (according to the Stone theorem \cite{ah}, $e_0$ is a cyclic
vector for $\cal L$), we have by the spectral theorem that
\begin{equation*}
s_k=\int_{-\infty}^{\infty} \g^k d(E_\g e_0,e_0)=({\cal L}^k
e_0,e_0)=(L^k)_{0,0},
\end{equation*}
so $s_k$ is the element of the matrix $L^k$ with indices $(0,0)$.

A key role in the inverse spectral problem method for Jacobi and
band operators is played by their Weyl functions (or Weyl matrices
for the band operators) \cite{yurm,osrjmp2}. In our case, it is
defined on the resolvent set of $\cal L$ as follows:
\begin{equation*}
\M=\M(z)=({\cal R}_z e_0,e_0),
\end{equation*}
where ${\cal R}=(z{\cal I}-{\cal L})^{-1}$ is the resolvent of
$\cal L$, $z \in \mathbb{C}$ and ${\cal I}$ is the identity
operator. It follows from the spectral theorem that
\begin{equation}
\label{mvj} \M(z)=\int_{-\infty}^{\infty} \frac{d\rho(\g)}{z-\g}.
\end{equation}
Thus, $\M(z)$ is the Stieltjes transform of
$\rho(\g)$\footnote{Similarly to \cite{nik}, we could call
$\rho(\g)$ the spectral function of $\cal L,$ but, as mentioned
above, the name ``spectral function" was previously reserved for
another object. The term ``spectral measure" is also ambiguous, so
we specify in \eqref{dro} what we mean by it.}. Using the Neumann
expansion for the resolvent, we find that the Weyl function of
$\cal L$ admits the following (formal, if the operator $\cal L$ is
unbounded) expansion at infinity
\begin{equation}
\label{neum} \M(z)=\sum_{k=0}^\infty \frac{s_k}{z^{k+1}}
\end{equation}
The sequence $S=\{s_k\}_{k=0}^{\infty}$ can be called the moment
sequence of the Weyl function of $\cal L$ \cite{osrjnd,osj}. Thus,
the inverse spectral problem for $\cal L$ admits the following
formalization: given $S,$ find the elements of $L$ (the operator
$\cal L$). Taking into account the above formula and using the
uniqueness of the Stieltjes transform, we may conclude that in the
context of our inverse spectral problem, the spectral measure $d
\rho(\g)$ and the Weyl function $\M(z)$ of $\cal L$ are two
equivalent objects; formula \eqref{solisp} provides the solution
(though for wider classes of band operators (i. e. operators
generated by infinite matrices containing a finite number of
nonzero diagonals), the Weyl function/matrix formalism \cite{yurm}
appears to be better than the use of $d\rho(\g) $ and its
generalizations, see e. g. \cite{osrjmp2}).
 Due to evenness of $d\rho(\g),$ we have
 $d\rho(\g)=0.5(d\rho(\g)+d\rho(-\g)),$ and \eqref{mvj} can be
 rewritten as follows:
\begin{equation}
\label{mveven}
 \M(z)=z\int_{-\infty}^{\infty}
 \frac{d\rho(\g)}{z^2-\g^2}.
\end{equation}
 Now turn back to the system \eqref{vl}-\eqref{seminf}. Assume that $(a_n(0))_{n=0}^{\infty}$ are given and our aim is
 to  find $(a_n(t))_{n=0}^{\infty}$ for all $t$. Since \eqref{vl}-\eqref{seminf} satisfies the
 Lax equation with the matrices $L(t)$ and $A(t),$
the equation
\begin{equation}
\label{lres} \overset \cdot R(t) = [R(t),A(t)],
\end{equation}
 is also hold, where the $R(t)$
 is a matrix representation of the resolvent ${\cal R} ={\cal
 R}(t)$ (we call $R(t)$ a resolvent matrix for $L(t)$).
Using this fact, we can find the evolution in time of the Weyl
function $m(z,t)$ of $L(t),$ see \cite{bert2} for the details
\begin{equation*}
\overset \cdot \M=\overset \cdot \M(z,t)=
\M(z,t)(z^2-a_0(t))-z=\M(z,t)(z^2-s_2(t))-z.
\end{equation*}

Let $d \rho(\g,t)$ be the spectral measure of ${\cal L}(t)$. Since
$\int_{-\infty}^{\infty}d\rho(\g,t)=1,$ the substitution of
$\displaystyle \M(z,t)=z\int_{-\infty}^{\infty}
 \frac{d\rho(\g,t)}{z^2-\g^2}$ into the above formula yields
\begin{equation*}
\overset \cdot \M(z,t) = z\int_{-\infty}^{\infty}
 \frac{(z^2-a_0(t))d\rho(\g,t)}{z^2-\g^2}-z\int_{-\infty}^{\infty}
 \frac{z^2-\g^2}{z^2-\g^2}d\rho(\g,t)=z\int_{-\infty}^{\infty}
 \frac{(\g^2-a_0(t))d\rho(\g,t)}{z^2-\g^2}.
\end{equation*}
Using the uniqueness of the Stieltjes transform, we get the
evolution in time of  $d\rho(\g,t):$
\begin{equation}
\label{eqrho} \overset \cdot
d\rho(\g,t)=(\g^2-a_0(t))d\rho(\g,t)=(\g^2-s_2(t))d\rho(\g,t).
\end{equation}
Thus we find
\begin{equation}
\label{drt}
d\rho(\g,t)=\frac{e^{\g^2t}d\rho(\g,0)}{e^{\int_{0}^t
a_0(\tau)d\tau}},
\end{equation}
from which, using again the normalization condition
$\int_{-\infty}^{\infty}d\rho(\g,t)=1,$ we finally get
\begin{equation}
\label{evrho} d\rho(\g,t)=\frac{e^{\g^2 t}d\rho(\g,0)}{
\int_{-\infty}^{\infty} e^{\g^2 t}d\rho(\g,0)}:=K(t)e^{\g^2
t}d\rho(\g,0).
\end{equation}
By our assumption, $L(0)$ is bounded, so $a(0)\in l_{\infty}$ and
the measure $d\rho(\g,0)$ has in fact a compact support. It
follows from \eqref{evrho} that this holds also for every $t\in
[0,T)$.

 For the moments $s_k=s_k(t)$, respectively, we obtain the
evolution equation
\begin{equation*}
\overset \cdot s_k(t) = s_{k+2}(t)-s_2(t)s_k(t), \quad k \in
\mathbb{Z_+},
\end{equation*}
(as we noted above, $s_k=s_k(t)=0$ for odd $k$) which has an
analytical solution \cite{osrjmp2}
\begin{equation}
\label{momev}
s_{2k}(t)=\frac{\displaystyle
\sum_{m=0}^{\infty}\frac{s_{2(k+m)}(0)t^m}{m!}}{\displaystyle
\sum_{m=0 }^{\infty}\frac{s_{2m}(0)t^m}{m!}}.
\end{equation}
Note that as mentioned in the Introduction, the theory of Jacobi
operators is related to the one of orthogonal polynomials (the
solution of inverse spectral problem for Jacobi operators via
orthogonal polynomials is contained in Chapter VII of \cite{ber}).
This relation in particular implies that that the measure
$d\rho(\g,0)$ is known for a number of initial data, for example,
the sequence
$(a_n(0))_{n=0}^{\infty}=(\sqrt{(n+1)/2})_{n=0}^{\infty}$
corresponds to $\displaystyle d\rho(\g)=
\frac{e^{-\g^2}}{\sqrt{\pi}}$ which is the orthogonality weight
for Hermite polynomials. Substituting $s_k(t)$ into
\eqref{solisp}, we get $a_n(t)$ which satisfy
\eqref{vl}-\eqref{seminf}.

We also note an important fact that the Weyl function of Jacobi
operator admits an expansion into infinite continued $J-$
fraction, see \cite{osj} and references thereafter, which in our
case takes the form
\begin{equation}
\label{jfr} \M(z,t) = \int_{-\infty}^{\infty}
\frac{d\rho(\g,t)}{z-\g} \sim \frac{1}{\displaystyle z -
\frac{a_0}{\displaystyle z-\frac{a_1}{\displaystyle \qquad
\ddots}}},
\end{equation}
This continued fraction approach was used by Nikishin in his study
of Volterra and Toda lattices \cite{nik,nisor}. The direct
spectral problem for Jacobi operators, which in our case amounts
to finding $d\rho(\g,0)$ (or $\M(z,0)$) from
$(a_n(0))_{n=0}^{\infty},\,$ may be interpreted as inversion of
the above continued fraction at $t=0,$ while getting the solution
$(a_n(t))_{n=0}^{\infty} $ of \eqref{vl}-\eqref{seminf} may be
regarded as expansion of $\M(z,t)$ into such $J-$ fraction.

Thus we briefly sketched the integration method for a
semi-infinite Volterra lattice by means of the inverse spectral
problem for Jacobi operators. Now consider the lattice \eqref{vl}
in the infinite case, i. e. when $n\in\mathbb{Z},$ assuming that
$\sup_n a_n(t) <\infty$ for all $t$.
 Here the Lax representation also holds with the matrices
$L(t)=(L_{ij}(t))_{i,j=-\infty}^{\infty}$ and
$A(t)=(A_{ij}(t))_{i,j=-\infty}^{\infty}$ having the same
structure as \eqref{lv} and \eqref{av}. Applying the duplication
method for Jacobi operators \cite{ber,zh,bg,ass}, one can
transform these infinite matrices $L(t)$ and $A(t)$ into
semi-infinite ones with matrix entries of size $2\times 2$ and
find that the system \eqref{lv} for $n \in \mathbb{Z}$ is
equivalent to the Lax equation
\begin{equation}
\label{lplus} \overset \cdot{\hat L}(t) = [\hat L(t), \hat A(t)],
\end{equation}
  where
\begin{eqnarray*}
\hat L(t)=(\hat L_{i,j}(t))_{i,j=0}^\infty=\begin{pmatrix} B_0& A_0& O&O&O&\dots\\
           A_0&O&A_1&O&O&\dots\\
           O&A_1&O&A_2&O&\dots\\
            O&O&\ddots&\ddots&\ddots&\dots
            \end{pmatrix},\\
            A_i=\begin{pmatrix} \sqrt{a_{-i-2}}& 0\\0
            &\sqrt{a_i}\end{pmatrix},\quad B_0=\begin{pmatrix}
            0&\sqrt{a_{-1}}\\\sqrt{a_{-1}}&0 \end{pmatrix},\quad
            O=\begin{pmatrix} 0&0\\0&0 \end{pmatrix};
\end{eqnarray*}
\begin{eqnarray*}
\hat A(t)=(\hat A_{i,j}(t))_{i,j=0}^\infty=\frac{1}{2}\begin{pmatrix}O& E_0& -C_0&O&O&O&\dots\\
           D_0&O&O&-C_1&O&O&\dots\\
           C_0&O&O&O&-C_2&O&\dots\\
            O&C_1&O&O&O&-C_3&\dots\\
            O&O&\ddots&\ddots&\ddots&\ddots&\ddots
            \end{pmatrix},\\
C_i=\begin{pmatrix} -\sqrt{a_{-i-3}a_{-i-2}}& 0\\0
            &\sqrt{a_ia_{i+1}}\end{pmatrix},\;D_0=\begin{pmatrix}
            0&-\sqrt{a_{-2}a_{-1}}\\\sqrt{a_{-1}a_{0}}&0
            \end{pmatrix},\; E_0=-D_0^{T}.
\end{eqnarray*}
For each $t,$ the matrix $ \hat L(t) $ generates the bounded (and
therefore selfadjoint) operator $ {\cal \hat L}(t) $ acting in
$l^2(\mathbb C^2, [0,\infty)),$  by which we denote the Hilbert
space of sequences $\{x_n\}_{n=0}^{\infty},\, \sum_{n=0}^{\infty}
(x_n, x_n) < \infty$. Also, let $\delta_0$ be the operator acting
from $\mathbb{C}^2$ to $l^2(\mathbb C^2, [0,\infty)),$ which takes
$x\in \mathbb{C}^2$ to the sequence $\{\delta_{0 n}x\}$. Then, if
$\{E_\g(t)\}$ is the spectral family of $ \hat {\cal L}(t), $ we
can define its spectral matrix \cite[Chapter VII]{ber} (spectral
measure)
\begin{equation}
\label{spmatr} d\hat \rho(\g,t)=
\begin{pmatrix} d\rho_{0,0}(\g,t) &
d\rho_{0,1}(\g,t)\\
d\rho_{1,0}(\g,t) & d\rho_{1,1}(\g,t)
\end{pmatrix},
\end{equation}
where $\hat\rho(\g,t)$ is a matrix representation of $\delta_0^*
E_\g(t)\delta_0$ (cf. with \eqref{dro}). The necessary as well as
sufficient conditions which the elements of \eqref{spmatr} must
satisfy in order to form the spectral matrix of $ \hat {\cal
L}(t)$ are contained in the Theorem 1 of \cite{zh}. The elements
of $\hat L(t)$ can be found from $ d\hat\rho(\g,t)$ through the
pseudo-orthogonalization procedure with respect to this measure
applied to the matrix polynomial system $\{\g^k
E\}_{k=0}^{\infty},\,$ where $E=(\delta_{i,j})_{i,j=0}^1,$ see
\cite{ber}.

Due to \eqref{lplus}, the equation similar to \eqref{lres} is hold
in this case, see e. g. \cite{osj} for description of a resolvent
matrix for $\hat {\cal L}(t)$. Using the approach contained in
\cite{bg} we find in the same manner as above that the spectral
matrix of $\hat{\cal L}(t)$ satisfies the evolution equation
\begin{eqnarray}
 \frac{d\hat\rho(\g,t)}{dt}=\frac{1}{2}\left(\sigma_3 \g^2 -
F(\g,t)\right)d\hat\rho(\g,0)\left(\sigma_3 \g^2 -
F^*(\g,t)\right)
\\
\text{where} \quad \sigma_3=\begin{pmatrix} -1 & 0\\0
            &1\end{pmatrix},\quad  F(\g,t)= \begin{pmatrix}
-a_{-2}(t)+a_{-1}(t)&-2\g \sqrt{a_{-1}(t)}\\2\g
\sqrt{a_{-1}(t)}&-a_{-1}(t)+a_{0}(t)
\end{pmatrix}
            \nonumber
\end{eqnarray}
We will call $d\hat\rho(\g,t)$ a spectral matrix  corresponding to
the infinite Volterra lattice. In what follows, we will also need
the Weyl matrix of $\hat{\cal L}(t):$
\begin{equation}
\label{wmp} \hat{\M}(z,t) = \begin{pmatrix} m_{0,0}(z,t) &
m_{0,1}(z,t)\\
m_{1,0}(z,t) & m_{1,1}(z,t)
\end{pmatrix}:= \begin{pmatrix} \displaystyle \frac{d\rho_{0,0}(\g,t)}{z-\g} &
\displaystyle \frac{d\rho_{0,1}(\g,t)}{z-\g}\\
\displaystyle\frac{d\rho_{1,0}(\g,t)}{z-\g} & \displaystyle
\frac{d\rho_{1,1}(\g,t)}{z-\g}
\end{pmatrix}.
\end{equation}

 Since $\sigma_3 $ and $F(\g,t) $ do
not commute with each other, we cannot obtain here the analogs of
formulas \eqref{drt}-\eqref{evrho}, so the above-considered
integration procedure for the semi-infinite Volterra lattice can
not be directly transferred to the infinite case. As shown in
\cite{bg}, the same situation holds for the classical Toda lattice
(see e. g. \cite{zh,os94} for another versions of the inverse
spectral problem method applied to the integration of the infinite
Toda and Volterra lattices).

The classical Volterra lattice has rich Hamiltonian structure.
Namely, define in the coordinates $(a_n)$ the quadratic Poisson
bracket $\{.,.\}_2$ having non-vanishing elements
\begin{equation*}
\{a_{n+1},a_n\}_2=a_{n+1}a_{n},
\end{equation*}
and, respectively, the cubic bracket $\{.,.\}_3$ with
non-vanishing elements
\begin{equation*}
\{a_{n+2},a_n\}_3=a_{n+2}a_{n+1}a_{n},\quad
\{a_{n+1},a_n\}_3=a_{n+1}a_n(a_{n+1}+a_n),
\end{equation*}
see \cite{dam,sur} for more details. Then it is easily verified
that the system \eqref{vl} can be written as follows:
\begin{equation*}
\overset \cdot a_n =\{H_0,a_n\}_3=\{H_1,a_n\}_2
\end{equation*}
with the Hamiltonians $H_0=\displaystyle \frac{1}{2}\sum_n \ln
a_n, $ and $H_1=\sum_n a_n$. Note that $H_1(t)=\displaystyle
\frac{\mathrm{Tr}(L^2(t))}{2},$ where $L(t)$ is the Lax matrix for
the system \eqref{vl} both in the semi-infinite and the infinite
cases (this formula holds as well for the finite lattices).
Moreover, this bi-Hamiltonian structure allows one to write all
the systems $(V_k)_{k=1}^\infty$ of Volterra hierarchy (Volterra
flow) starting from $V_1,$ which is the original Volterra lattice
\eqref{vl} in the following manner  (see \cite{osrjnd} and
references thereafter):
\begin{eqnarray}
\label{2ham}
V_k:=\quad \overset \cdot a_n =\{H_{k-1},a_n\}_3=\{H_k,a_n\}_2,\\
\text{where} \quad H_k(t)=\displaystyle
\frac{\mathrm{Tr}(L^{2k}(t))}{2}. \nonumber
\end{eqnarray}
Faddeev and Takhtajan in \cite{ft1} found another ``interesting"
\cite{ft2} structure for \eqref{vl} generated by the following
Poisson bracket:
\begin{equation}
\label{ft} \{a_{n+1},a_n\}_{ft} = a_{n+1}a_n(a_{n+1}+a_n-4),\quad
\{a_{n+2},a_n\}_{ft} = a_{n+2}a_{n+1}a_n.
\end{equation}
As mentioned in \cite{fv}, it may be regarded as a lattice
counterpart of the Magri bracket. It defines a lattice Virasoro
algebra known as Faddeev-Takhtajan-Volkov algebra \cite{hik1,fv}.
Using \eqref{ft}, we get the following representation for
\eqref{vl}
\begin{equation}
\label{vft} \overset \cdot a_n =\{H_0,a_n\}_{ft}.
\end{equation}
Formulas \eqref{2ham}, \eqref{vft} are valid in the finite,
semi-infinite and infinite cases.
 As we see, the logarithmic Hamiltonian $H_0$ is present
in both \eqref{2ham} and \eqref{vft}. Our further goal will be to
find the conditions for its finiteness in terms of the
above-considered spectral measures corresponding to the lattices.
%%%%%%%%%%%%%%%%%%%%%%%%%%%%%%%%%%%%%%%%%%%%%%%%%%%%%
%%%%%%%%%%%%%%%%%%%%%%%%%%%%%%%%%%%%%%%%%%%%%%%%%%%%%
\section{Orthogonal polynomials. The Nikishin theorem} We proceed
with a brief discussion of results of previously mentioned work
\cite{nik}, essential for our subject. They are mostly related to
improvement of the classical results about the asymptotic
properties of orthogonal polynomials on the unit circle and the
unit interval \cite{seg,ger} (Szeg\H{o} asymptotic formulas)
aiming to utilize them for the inverse scattering method for
Jacobi operators. One of the tasks here was to study the cases in
which the spectral measures of the latter satisfy the Szeg\H{o}
condition.

According to \cite{nik}, first consider the operator $\cal L$ such
that its spectrum entirely consists of an infinite number of
points in the interval $[-1,1]$ (e. g. the operator generated by
\eqref{lv} with $a_n=1/4$ fits this description \cite{nisor}). Let
$d\rho(\g)$ be the corresponding measure \eqref{dro}. Assuming
that $\rho(-1)=0,$ we define a measure $d \sigma(\theta)$ on
$[0,2\pi]$ in the following way:
\begin{equation}
\label{sig} \sigma(\theta)=\begin{cases} -\rho(cos(\theta)),\quad 0\le \theta \le \pi;\\
\;\;\rho(cos(\theta)),\quad \pi\le \theta \le 2\pi.
\end{cases}
\end{equation}
As we see, $\sigma(\theta)=-\sigma(2\pi-\theta),\;  0\le \theta
\le \pi$. In \cite{nik} the measure $d\sigma(\theta)$ constructed
according to \eqref{sig} is called symmetric. As mentioned above,
in our case the measure $d\rho(\g)$ is even, and additionally we
have $\sigma(\theta)+\sigma(\pi-\theta)=2\sigma(\pi/2),\; 0\le
\theta \le
\pi,\;\sigma(\theta)+\sigma(2\pi-\theta)=2\sigma(3\pi/2),\quad
\pi\le \theta \le 2\pi\,$.

Let $\{\Phi_n(z)\}_{n=0}^{\infty}$ be the polynomials, $\deg
\Phi_n(z)=n,$ which are orthogonal on the unite circle with
respect to $d\sigma(\theta)$ and have the leading coefficient one.
Denote by
$c_{k}=\frac{1}{2\pi}\int_{0}^{2\pi}e^{ik\theta}d\sigma(\theta),\;
k \in \mathbb{Z}$ the trigonometric moments of $d\sigma(\theta)$.
Since the latter is symmetric, $c_k \in \mathbb{R}$ for all $k$.
Also, the fact that $d\sigma(\theta)$ is even for $ 0\le \theta
\le \pi$ and $\pi\le \theta \le 2\pi$ implies that $c_k=0$ for odd
$k$.

The polynomials $\Phi_n(z)$ can be calculated as follows, see e.
g. \cite[Chapter 5] {ah}
\begin{equation*}
\Phi_n(z)=\frac{1}{T_{n-1}}\left\vert\begin{array}{cccc} c_0&
c_1&\ldots &c_n\\
c_{-1}&c_0&\ldots &c_{n-1}\\
\vdots&\vdots&\ldots &\vdots\\
c_{-n+1}&c_{-n+2}&\ldots &c_{1}\\
1&z&\ldots &z^n
\end{array}\right\vert,
\end{equation*}
where $T_n$ is a Toeplitz determinant $|c_{i-j}|_0^n$. It is
easily verified that
\begin{equation*}
\frac{1}{2\pi}\int_{0}^{2\pi}\Phi_n(e^{i\theta})
\overline{\Phi_m(e^{i\theta})}
d\theta=\frac{1}{2\pi}\int_{0}^{2\pi}\Phi_n(e^{i\theta})
\Phi_m(e^{-i\theta}) d\theta=\frac{T_n}{T_{n-1}}\delta_{n,m}.
\end{equation*}
The following recurrence relations are fulfilled for these
polynomials (Szeg\H{o}'s recurrence):
\begin{eqnarray}
\label{srec}
\Phi_{n+1}(z)=z \Phi_n(z)-\alpha_n \Phi_n^*(z),\quad\\
 \text{where}\quad \Phi_n^*(z)=z^n
\Phi_n\left(\frac{1}{z}\right),\quad \alpha_n = -\Phi_{n+1}(0)\in
\mathbb{R}. \nonumber
\end{eqnarray}
In \cite{nik}, the numbers $(\alpha_n)$ were called the circular
parameters (of the operator $\cal L$) though later on the names
``reflection coefficients" and ``Verblunsky coefficients"
\cite{sim1,simop} came into usage. In our case, since $d\rho(\g)$
and, respectively, $d\sigma(\theta )$ have infinitely many growth
points, we have that $T_n>0$ for all $n,$ (the moment sequence
$\{c_k\}_{k=-\infty}^{\infty}$ is positive) which in turn implies
that $0\le |\alpha_n| <1$, see \cite[Chapter V]{ah} for more
details. Using the fact that in our case the odd moments $c_k$ are
zero, one can check that $\alpha_n=0$ for even $n$. The
coefficients of $L$ can be expressed in terms of $(\alpha_n)$ as
follows \cite{seg,nik}:
\begin{equation}
\label{aalpha}
a_n=\frac{(1-\alpha_{n-1})(1-\alpha_n^2)(1+\alpha_{n+1})}{4},\; n
\in \mathbb{N},\quad
a_0=\frac{(1-\alpha_0^2)(1+\alpha_1)}{2}=\frac{(1+\alpha_1)}{2}.
\end{equation}
As shown by Geronimus, see \cite{ah} and references thereafter,
for the measure $d\tilde \sigma(\theta)$ which solves the
trigonometric moment problem, the Szeg\H{o} condition
\begin{equation}
\label{segtrig} \int_{0}^{2\pi} \ln \tilde \sigma
'(\theta)d\theta
>-\infty,
\end{equation}
(where $\tilde \sigma '(\theta)$ is the derivative of the
absolutely continuous part of $\tilde \sigma (\theta)$)
 is equivalent to
\begin{equation*}
\prod_{n=1}^{\infty}(1-|\tilde \alpha_n|^2)<\infty,
\end{equation*}
which can be rewritten as
\begin{equation}
\label{alcon}
\sum_{n=1}^{\infty}|\tilde \alpha_n|^2<\infty,
\end{equation}
where $\tilde \alpha_n $ are the Verblunsky coefficients of the
corresponding Szeg\H{o}'s recurrence. Turning back to our case, we
find using \eqref{aalpha} that in terms of the operator $\cal L,$
the Szeg\H{o} condition can be written as follows
\begin{equation}
\label{seg1}  \int_{-1}^1 \frac{\ln \rho'(\g)}{\sqrt{1-\g^2}}d\g
<-\infty,
\end{equation}
while \eqref{alcon} is equivalent to
\begin{equation}
\label{prodcond} 0<\prod_{n=0}^{\infty} 4a_n < \infty,
\end{equation}
see also \cite[Supplement]{seg}.

 In \cite{nik}, Nikishin considered the situation, when in
addition to the points inside $[-1,1],$ the spectrum of $\cal L$
consists of a finite number of points (eigenvalues) located
outside this interval. He obtained a generalization of Szeg\H{o}
asymptotic formulas to this case, and this result plays an
important role in the inverse scattering method for Jacobi
operators proposed in this paper. As a matter of fact, he studied
more general case of selfadjoint operators ${\cal L}^\dagger$
generated by Jacobi matrices $L^\dagger$ (such operators appear e.
g. in the Lax pair which corresponds to Toda lattice):
\begin{equation}
\label{hatl}
L^\dagger = \begin{pmatrix} v_0& w_0& 0&0&0&\dots\\
           w_0&v_1&w_1&0&0&\dots\\
           0&w_1&v_2&w_2&0&\dots\\
            0&0&\ddots&\ddots&\ddots&\dots
            \end{pmatrix},
\end{equation}
where $v_n$ and $w_n$ are real. In the context of inverse
scattering problem for discrete Shturm-Liouville (Jacobi)
operators, the eigenvalues laying outside $[-1,1]$ correspond to
bound states which, in particular, are related to solitons in the
lattices \cite{to}. An important result established in \cite{nik}
(Theorem 2) is that for these operators ${\cal L}^\dagger$ there
exists a number $M\ge 1$ such that the spectrum of operator ${\cal
L}_M^\dagger$ generated by the ``shifted" matrix
\begin{equation*}
L_M^\dagger = \begin{pmatrix} v_M& w_M& 0&0&0&\dots\\
           w_M&v_{M+1}&w_{M+1}&0&0&\dots\\
           0 &w_{M+1}&v_{M+2}&w_{M+2}&0&\dots\\
            0&0&\ddots&\ddots&\ddots&\dots
            \end{pmatrix}.
\end{equation*}
has no points of its spectrum (eigenvalues) outside $[-1,1]$. By
taking the spectral measures \eqref{dro} $ d\mu^\dagger(\g)$ and $
d\mu_M^\dagger(\g)$ of the operators ${\cal L}^\dagger$ and ${\cal
L}_M^\dagger$ respectively, Nikishin considered the expansion of
the Weyl function $\M^\dagger (z):=\int_I \frac{
d\mu^\dagger(\g)}{z-\g}$ of ${\cal L}^\dagger,$ where $I$ is the
support of its spectrum, into the following $J-$ fraction:
\begin{equation*}
\M^\dagger(z) = \frac{1}{\displaystyle z - v_0 -
\frac{w_0^2}{\displaystyle z-v_1-\frac{w_1^2}{\displaystyle \qquad
\;\frac{\ddots}{\quad z - \;v_{M-1} -\; w_{M-1}^2
\M_M^\dagger(z)}}}},
\end{equation*}
where $\M_M^\dagger(z)=\int_{-1}^{1}
\frac{d\mu_M^\dagger(\g)}{z-\g}$ is the Weyl function of ${\cal
L}_M^\dagger $. As well known, see e. g. \cite{ah}, the
convergents of this continued fraction can be written in terms of
polynomials of the second and the first kind associated with
Jacobi matrix $L^\dagger$ (as their ratios). Using this fact as
well as some results from the complex analysis (mainly related to
boundary values of analytic functions) it was shown that the
following relation is hold
\begin{equation*}
\ln(\mu^{\dagger}(\g))'=\ln (\mu^{\dagger}_M(\g))'+f(\g),
\end{equation*}
where $f(\g)$ is a certain function  such that
$\frac{f(\g)}{\sqrt{1-\g^2}}\in L_1[-1,1]$. From this it
immediately follows that the conditions
\begin{equation}
\label{seghat} \int_{-1}^1 \frac{\ln
(\mu^\dagger(\g))'}{\sqrt{1-\g^2}}\,d\g>-\infty \quad
\end{equation}
and
\begin{equation*}
\int_{-1}^1 \frac{\ln
(\mu^\dagger_M(\g))'}{\sqrt{1-\g^2}}\,d\g>-\infty
\end{equation*}
are equivalent. The latter one, in turn, is equivalent to
\begin{equation*}
 0<\prod_{n=M}^{\infty} 4w_n^2 < \infty,
\end{equation*}
cf. with \eqref{seg1}-\eqref{prodcond}. Therefore, since the
convergence of an infinite product is not affected by inserting a
finite number of factors, the following result is hold:
%%%%%%%%%%%%%%%%%%%%%%%%%%%%%%%%%%%%%%%%%%%%%%%%%%%%%%%%%%%%%%%%
\begin{thm*}
Assume that the spectrum of the operator ${\cal L}^\dagger$
generated by Jacobi matrix \eqref{hatl} consists of an infinite
number of points in the interval [-1,1] and a finite number of
points (eigenvalues) outside this interval. Then the corresponding
measure $d\mu^\dagger(\g),$ defined as in \eqref{dro}, satisfies
the Szeg\H{o} condition \eqref{seghat} if and only if
\begin{equation}
\label{nikprod}
 0<\prod_{n=0}^{\infty} 4w_n^2 < \infty,
\end{equation}.
\end{thm*}
%%%%%%%%%%%%%%%%%%%%%%%%%%%%%%%%%%%%%%%%%%%%%%%%%%%%%%%%%%%%%%%%
This formulation is the same as in \cite{nik} though some remarks
are needed with regards to the spectrum of ${\cal L}^\dagger$. Its
infiniteness on $[-1,1]$ was used in the proof of the fact that
for certain $M$ the spectrum of ${\cal L}_M^\dagger$ lies entirely
within this interval. At the same time, the Szeg\H{o} condition
implies that ${\cal L}^\dagger$ must have additional spectral
properties. Obviously, the derivative $\mu^\dagger(\g)\,'$ must
exist and be positive almost everywhere on $[-1,1]$. We can
conclude then, that the essential spectrum of ${\cal L}^\dagger$
is $[-1,1],$ meanwhile, we cannot guarantee that its singularly
continuous spectrum is empty. As also mentioned in \cite{nik}, the
latter condition is essential in the context of the inverse
scattering method for such operators, because if it is fulfilled
(i. e. the spectrum on $[-1,1]$ is purely absolutely continuous),
then the operator ${\cal L}^\dagger$ can be uniquely reconstructed
from the scattering data . It was shown in \cite{gus1} that this
condition holds if
\begin{equation}
\label{marc} \sum_{n=0}^\infty n
\left(\left|w_n^2-\frac{1}{4}\right| + |v_n|\right)<\infty.
\end{equation}
In some sources, see e. g. \cite{aptnik}, by the analogy with the
Sturm-Liouville differential operator case, \eqref{marc} is called
the Marchenko condition. It also ensures the at most finiteness of
the point spectrum of  ${\cal L}^\dagger$ (the eigenvalues located
outside $[-1,1]$), see also \cite{gk}, Theorem 1. In \cite{nik} it
was shown that the Marchenko condition can be weakened to some
extent to ensure both the pure absolute continuity of the spectrum
of ${\cal L}^\dagger$ on $[-1,1]$ and the finiteness of its point
spectrum. Another result contained in \cite{nik} is that if, in
addition to \eqref{prodcond}, $\sum_n v_n^2 =\infty,$ then the
corresponding operator ${\cal L}^\dagger$ has an infinite number
of eigenvalues outside $[-1,1]$. It should be noted that since the
publication of \cite{nik}, the theory of Jacobi operators, where
the influence of asymptotic properties of their coefficients on
the structure of their spectra is among the main issues, have been
a subject of active studies, some details will be mentioned below;
see \cite{sim1}-\cite{simop} for an intermediate summary of
results.
%%%%%%%%%%%%%%%%%%%%%%%%%%%%%%%%%%%%%%%%%%%%%%%%%%%%%%
%%%%%%%%%%%%%%%%%%%%%%%%%%%%%%%%%%%%%%%%%%%%%%%%%%%%%%%%%%%%%%%
\section{Description of logarithmic Hamiltonians} Now consider
the operator ${\cal L}$ generated by \eqref{hatl} with the zero
main diagonal, such that its spectral measure $d \rho(\g)$ (which
in this case is even) matches the conditions of Nikishin's
theorem. Then, because of the multiples 4 in \eqref{nikprod}, we
find that $\displaystyle \sum_{n=0}^\infty \ln w_n^2 = -\infty$ .
At the same time, if we set
\begin{equation}
\label{transmes} d\tilde \rho(\xi)= d\rho(2\g).
\end{equation}
then the Szeg\H{o} condition \eqref{seghat} for $d\hat\rho(\g)$
turns into
\begin{equation*}
\label{seg2} \int_{-2}^2 \frac{\ln \tilde\rho\,
'(\xi)}{\sqrt{4-\xi^2}}\,d\xi>-\infty
\end{equation*}
The moments $\tilde s_k := \int_{-\infty}^{\infty}\xi^k d\tilde
\rho(\xi)$ and $s_k = \int_{-\infty}^{\infty}\g^k d\rho(\g)$ are
then related to each other as $\tilde s_k = 2^{k} s_k, \, (\tilde
s_0 = s_0=1)$. For the determinants $\tilde \Delta_n := |\tilde
s_{i+j}|_{i,j=0}^{n}$ and $\Delta_n := | s_{i+j}|_{i,j=0}^{n}$ we
have $\tilde \Delta_n = 2^{(n+1)n}\Delta_n>0$ and, respectively,
\begin{equation}
\label{tila} \tilde
w_n^2=\frac{\tilde\Delta_{n+1}\tilde\Delta_{n-1}}{\tilde\Delta^2_n}=4\frac{\Delta_{n+1}\Delta_{n-1}}{\Delta^2_n}=4w_n^2,
\end{equation}
so the condition \eqref{nikprod} turns into
\begin{equation*}
 0<\prod_{n=0}^{\infty} \tilde w_n^2 < \infty.
\end{equation*}.

 Now take an even probability measure $d\tilde \rho(\xi)=d\tilde
\rho(\xi,0)$ which, as in Nikishin's theorem, has infinitely many
growth points in the interval $[-2,2],$ a finite number of masses
located on $\mathbb{R}\backslash [-2,2]$ and satisfies the
Szeg\H{o} condition
\begin{equation}
\label{seg20}
 \int_{-2}^2 \frac{\ln \tilde\rho\,'(\xi)}{\sqrt{4-\xi^2}}\,d\xi>-\infty
\end{equation}
 on $[-2,2]$.
 Then, by the Hamburger
theorem, the moment sequence $(\tilde s_k(0))_{k=0}^{\infty},$ of
$d\tilde \rho(\xi)$ is positive, so we can construct the sparse
Jacobi matrix
\begin{equation*}
\tilde L(0) = \begin{pmatrix} 0& \sqrt{\tilde a_0(0)}& 0&0&0&\dots\\
           \sqrt{\tilde a_0(0)}&0&\sqrt{\tilde a_1(0)}&0&0&\dots\\
           0&\sqrt{\tilde a_1(0)}&0&\sqrt{\tilde a_2(0)}&0&\dots\\
            0&0&\ddots&\ddots&\ddots&\dots
            \end{pmatrix},
\end{equation*}
where $\displaystyle \tilde a_n(0) =
            \frac{\tilde\Delta_{n+1}(0)\tilde\Delta_{n-1}(0)}{\tilde\Delta^2_n(0)},$
such that $d\tilde \rho(\xi,0)$ is the spectral measure of the
corresponding operator $\tilde {\cal L}(0).$ Acting similarly as
above, we find that
\begin{equation*}
 0<\prod_{n=0}^{\infty} \tilde a_n(0) < \infty.
\end{equation*}

 Next, define the measure $d\tilde \rho(\xi,t)$ as in
\eqref{evrho}: $d\tilde \rho(\xi,t) = K(t)e^{\xi^2 t}d\tilde
\rho(\xi,0).$ It follows from this definition that such measure
has the same essential support $[-2,2]$ and the mass points
outside this interval as $d\tilde \rho(\xi,0) $ does, and we also
get
\begin{equation*}
\int_{-2}^2 \frac{\ln \tilde\rho\,'(\xi,t)}{\sqrt{4-\xi^2}}\,d\xi
= 2\pi t + \int_{-2}^2 \frac{\ln
\tilde\rho\,'(\xi,0)}{\sqrt{4-\xi^2}}\,d\xi + \pi \ln K(t)>-\infty
\end{equation*}
Hence
\begin{equation}
\label{tilprod}
 0<\prod_{n=0}^{\infty} \tilde a_n(t) < \infty
\end{equation}
for $t>0$ and we get the result which can be regarded as a
corollary of the Nikishin theorem.
%%%%%%%%%%%%%%%%%%%%%%%%%%%%%%%%%%%%%%%%%%%%%%%%%
\begin{theorem}
\label{th1}
Each even probability measure $d\tilde \rho(\xi)$ with
essential support coinciding with the interval $[-2,2],$ having a
finite number of mass points on $\mathbb{R}\backslash [-2,2]$ and
satisfying \eqref{seg20}, generates the semi-infinite Volterra
lattice
\begin{equation}
\label{tilv}
\overset \cdot {\tilde a}_n = \tilde a_n (\tilde
a_{n+1}-\tilde a_{n-1}), \qquad n \in \mathbb Z_+,
\end{equation}
with a finite logarithmic Hamiltonian $\tilde H_0 =\sum_n \ln
\tilde a_n $ by means of the following procedure
\begin{itemize}
\item[a)] Define the measure $d\tilde \rho(\xi,t):$
\begin{equation}
\label{rhohi} d\tilde\rho(\xi,t)=\frac{e^{\xi^2 t}d\tilde
\rho(\xi)}{ \int_{-\infty}^{\infty} e^{\xi^2 t}d\tilde \rho(\xi)}.
\end{equation}
\item[b)] Calculate the moments of $d\tilde\rho(\xi,t)$: $ \tilde
s_k(t)=\int_{-\infty}^{\infty} \xi^k d\tilde\rho(\xi,t), k \in
\mathbb{Z}_+, \; (\tilde s_{2k}=0).$ \item[c)] Find $\tilde
a_n(t):$
\begin{equation*}
\tilde a_n(t) =
            \frac{\tilde\Delta_{n+1}(t)\tilde\Delta_{n-1}(t)}{\tilde\Delta^2_n(t)},
            \quad\text{where}\; \tilde \Delta_n :=
|\tilde s_{i+j}(t)|_{i,j=0}^{n}.
\end{equation*}
\end{itemize}
\end{theorem}
%%%%%%%%%%%%%%%%%%%%%%%%%%%%%%%%%%%%%%%%%%%%%%%%%%%%%%%%%%%%%%%%%%%%%%
Now consider the general case of \eqref{vl}-\eqref{seminf} with a
finite logarithmic Hamiltonian $H_0$ . Its finiteness obviously
implies that $a_n(t) \to 1$ as $t \to \infty$, so the operator
${\cal L}(t)$ can be regarded as a compact perturbation of the
operator ${\cal J}_0$ generated by the ``free" matrix
\begin{equation*}
 J_0 = \begin{pmatrix} 0& 1& 0&\dots\\
          1&0&1&\dots\\
            0&\ddots&\ddots&\ddots
            \end{pmatrix},
\end{equation*}
As known, see e. g. \cite{nisor}, the spectrum of ${\cal J}_0$ is
$[-2,2]$ with the corresponding measure $d\rho_0(\g)=\displaystyle
\frac{1}{2\pi} \sqrt{4-\g^2}$ is purely absolutely continuous on
this interval. Then, according to the Weyl theorem, the spectrum
of ${\cal L}={\cal L}(t)$ is $[-2,2]$ plus no more than a
countable number of points (eigenvalues) outside $[-2,2]$ with a
point of accumulation only at the endpoints of this interval.
Using these findings and applying the previous results concerning
the inverse spectral problem for Jacobi operators, in particular
\eqref{solisp}, we find a description of the systems
\eqref{vl}-\eqref{seminf} having a finite logarithmic Hamiltonian
in terms of the measure $d\rho(\g,t)$ corresponding to the
operator ${\cal L}$ and its moments.
%%%%%%%%%%%%%%%%%%%%%%%%%%%%%%%%%%%%%%%%%%%%%%%%%%
\begin{theorem}
\label{th2} There is a one-to-one correspondence between the
systems \eqref{vl}-\eqref{seminf} with a finite logarithmic
Hamiltonian and the even probability measures $d\rho(\g,t)$ having
the essential support $[-2,2]$ and at most countable number of
mass points with accumulation points $\{-2,2\}$ such that for its
moments $s_k=s_k(t)=\int_{-\infty}^{\infty} \g^k d\rho(\g,t), k
\in \mathbb{Z}_+, \; (s_{2k}=0) $ the following condition is hold:
\begin{equation}
\label{delcond} 0<\lim_{n\to \infty}
\frac{\Delta_{n+1}}{\Delta_{n}}<\infty, \quad\text{where}\quad
\Delta_n=\det(s_{i+j})_{i,j=0}^n.
\end{equation}
Each of these systems has the Lax operator ${\cal L}(t)$ with the
spectral measure $d\rho(\g,t),$  possessing such properties.
Conversely, taking the measure $d\rho(\g,t)$ and applying
\eqref{solisp}, one can construct the system
\eqref{vl}-\eqref{seminf} having a finite logarithmic Hamiltonian.
\end{theorem}
%%%%%%%%%%%%%%%%%%%%%%%%%%%%%%%%%%%%%%%%%%%%%%%%%%%%%%%%%
It follows from \eqref{evrho}, that if the essential support of
$d\rho(\g,t)$ for $t=0$ is $[-2,2]$ then so is for each $t\in
[0,T)$ and its mass points are invariant with respect to $t$. Note
that the latter also follows from the fact that for the system
\eqref{vl}-\eqref{seminf}, corresponding to $d\rho(\g,t),$ the
matrix $A(t)$ \eqref{av} is skew-symmetric, which implies that the
eigenvalues of $L(t)$ defined by \eqref{lv} are time-invariant.
This is a general fact in the theory of integrable systems, see e.
g. \cite{bogo,to}. At the same time, \eqref{evrho} does not
automatically imply that once \eqref{delcond} is fulfilled for
$t=0,$  it holds as well for each $t\in [0,T),$ and a further
study is needed here.

Now, note that the Nikishin theorem does not cover the case when
the point spectrum of ${\cal L}^\dagger$ is infinite. As
previously mentioned, one of the results of \cite{nik} was
generalization of Szeg\H{o} asymptotic formulas to the case then
${\cal L}^\dagger$ has a finite point spectrum located outside
$[-1,1]$. In \cite{pyd}, Peherstorfer and Yuditskii showed that
under some additional condition this result can be generalized to
the case of infinite point spectrum. More precisely, they
considered the bounded Jacobi operator $\cal J$ generated by
Jacobi matrix
\begin{equation*}
 J = \begin{pmatrix} q_0& p_1& 0&\dots\\
          p_1&q_1&p_2&\dots\\
            0&\ddots&\ddots&\ddots
            \end{pmatrix},
\end{equation*}
having $[-2,2]$ as its essential spectrum and an infinite set of
eigenvalues
$\{\g_j^{+}\}_{j=1}^{\infty},\,\{\g_j^{-}\}_{j=1}^{\infty}$ such
that
\begin{equation}
\label{eigpm}
\g_{1}^{-}<\g_{2}^- <\dots < -2 < 2 < \dots <
\g_{2}^{+} < \g_{1}^{+}.
\end{equation}
They found that if the latter satisfy the condition
\begin{equation}
\label{pyds} \sum_{j=1}^{\infty} \sqrt{(\g_j^-)^2-4} +
\sqrt{(\g_j^+)^2-4} < \infty
\end{equation} and the Szeg\H{o} condition is hold for its spectral measure
$d \rho(\g, {\cal J})$ on $[-2, 2],$ then the asymptotic formulas
of Szeg\H{o} type are hold as well for the corresponding system of
orthogonal polynomials. In their paper, \eqref{pyds} was called
``Blaschke condition" (after the transformation from the real line
to the unit circle, used in \cite{pyd}, \eqref{pyds} becomes the
standard Blaschke convergence condition), but this name did not
take root among the specialists in the subject, so we call it
Peherstorfer-Yuditskii condition. It was also found that if the
Szeg\H{o} condition holds, while the one given by \eqref{pyds}
fails, then $\sum_{n=1}^{\infty} \ln p_n = \infty$ (in \cite{pyd},
it is written in the remark to the first Lemma of this paper).
These results were taken into consideration by Simon and
Zlato\v{c} in \cite{sz}, where it was found that if both Szeg\H{o}
and Peherstorfer-Yuditskii conditions are fulfilled for $d
\rho(\g, {\cal J}),$  then $\sum_{n=1}^{\infty}\ln p_n$ is finite
and ${\cal J}-{\cal J}_0$ is in Hilbert-Schmidt class, that is
\begin{equation*}
\sum_{n=1}^{\infty} (p_n-1)^2 + q_{n-1}^2 < \infty,
\end{equation*}
and, additionally, $\sum_{n=0}^{\infty} q_n $ is convergent (see
also Theorem 13.8.9 of \cite{simop}).

On the other hand, as shown by Hundeltmark and Simon \cite{hunds},
if the essential spectrum of ${\cal J}$ is $[-2,2],$ and its
eigenvalues satisfy \eqref{eigpm}, then
\begin{equation*}
\sum_{j=1}^{\infty} \sqrt{(\g_j^-)^2-4} + \sqrt{(\g_j^+)^2-4} <
4\sum_{n=1}^{\infty}|p_n-1| + |q_{n-1}|,
\end{equation*}
(in a recent work by Laptev, Loss and Schimmer \cite{lap}, this
estimate was improved by replacing $|p_n-1|$ with $\max\{p_n-1,
0\}$). Therefore, if we take the operator $\cal J,$ such that it
is in the trace class:
\begin{equation*}
\sum_{n=1}^{\infty} |p_n-1| + |q_{n-1}| < \infty,
\end{equation*}
which obviously implies that $ \sum_{n=1}^{\infty} \ln p_n $ is
finite, then the Peherstorfer-Yuditskii condition is fulfilled for
its eigenvalues. Also, as shown by Killip and Simon, see
\cite[Theorem 2]{kils}, if $\cal J$ is in the trace class, then
its spectral measure satisfies the Szeg\H{o} condition on the
interval $[2,2]$ thus confirming the Nevai's conjecture for Jacobi
matrices. Hence, both of the conditions are valid for such
operators $\cal J.$

All these results are applicable to the above operator ${\cal L}$,
where we take into account that in this situation its spectral
measure is even, so $\xi_{j}^{-}=-\xi_{j}^{+}$. Therefore, if we
take an even probability measure $d\rho(\xi)$ having the same
properties as $\tilde d \rho(\xi)$ except that its point spectrum
is infinite and satisfies the Peherstorfer-Yuditskii condition,
and acting similarly as in the case of $d \tilde\rho(\xi),$ we get
the following extension of the Theorem \ref{th1}.
%%%%%%%%%%%%%%%%%%%%%%%%%%%%%%%%%%%%%%%%%%%%%%%%%%%%%%%%%%%%%%%
\begin{theorem}
\label{th3} Each even probability measure $d\rho(\xi)$ with
essential support coinciding with the interval $[-2,2],$
satisfying there the the Szeg\H{o} condition,  having, possibly
infinite, point spectrum $\{\xi_j^{+}\},\,\{\xi_j^{-}\}, \;
\xi_j^{-} = - \xi_j^{+} $ satisfying the Peherstorfer-Yuditskii
condition
\begin{equation}
\label{pu}
\sum_{j}\sqrt{(\xi_j^+)^2-4} < \infty,
\end{equation}
generates the semi-infinite Volterra lattice with a finite
logarithmic Hamiltonian and such that
\begin{equation*}
\sum_{n=0}^{\infty} (a_n(t)-1)^2 < \infty,
\end{equation*}
 by means of the procedure specified in
Theorem \ref{th1} with $d\rho(\xi)$ instead of $d\tilde
\rho(\xi)$.

Conversely, for each system \eqref{vl}-\eqref{seminf}, such that
the series $\displaystyle \sum _{n=0}^{\infty}(a_n(t) - 1)$
absolutely converges for each $t$, the spectral measure of the
corresponding Lax operator ${\cal L}(0)$ has the same properties
as the measure $d\rho(\xi)$.
\end{theorem}
%%%%%%%%%%%%%%%%%%%%%%%%%%%%%%%%%%%%%%%%%%%%%%%%%%%%%%%%%%%%%%%%%
Using this result, we can make the following conclusion with
regards to the above condition \eqref{delcond}.
\begin{corollary}
If the measure $d\rho(\g) $ satisfies the conditions of Theorem
\ref{th3}, then for the moments $s_k(t) = \int_{-\infty}^{\infty}
\g^k d\rho(\g,t) $ of the measure $d\rho(\g,t)$ defined by
\eqref{evrho}, the condition \eqref{delcond} of Theorem \ref{th2}
is fulfilled for $t\ge 0$.
\end{corollary}
\begin{corollary}
Each even probability measure $d\rho_c(\xi)$ with essential
support coinciding with the interval $[-c,c],\; c>0$ satisfying
there the Szeg\H{o} condition
\begin{equation}
\label{segc} \int_{-c}^c \frac{\ln
\rho_c\,'(\xi)}{\sqrt{c^2-\xi^2}}\,d\xi
> -\infty,
\end{equation}
having, possibly infinite, point spectrum
$\{\xi_j^{+}\},\,\{\xi_j^{-}\}, \; \xi_j^{-} = - \xi_j^{+} $ such
that
\begin{equation}
\label{yudc} \sum_{j}\sqrt{(\xi_j^+)^2-c^2} < \infty,
\end{equation}
generates the semi-infinite lattice \eqref{vl} having the finite
regularized logarithmic Hamiltonian
\begin{equation*}
H_0^r := \sum_{n=0}^{\infty} (\ln (a_n) - r(c)),
\quad\text{where}\quad r(c)=2\ln \left(\frac{c}{2}\right).
\end{equation*}
\end{corollary}
Indeed, using the change of variable $\g = 2\xi/c$ in $d
\rho_c(\xi)$ and applying Nikishin's theorem, we find that
\eqref{segc} is equivalent to
\begin{equation*}
 0<\prod_{n=0}^{\infty} \frac{4}{c^2}a_n(0) < \infty
\end{equation*}
where $a_n(0)$ are defined according to \eqref{solisp}, where
$\Delta_n =|s_{i+j}|_{i,j=0}^{n} $ with
$s_{i+j}=\int_{-\infty}^{\infty} \xi^{i+j} d\rho_c(\xi).$ Next,
for $t\ge 0,$ define the measure
$d\rho_c(\xi,t)=K(t)\exp(\xi^2t)d\rho_c(\xi),$ with the
normalizing factor
$K(t)=(\int_{-\infty}^{\infty}\exp(\xi^2t)d\rho_c(\xi))^{-1}$.
Similarly as above, we check that \eqref{segc} is fulfilled for
$d\rho_c(\xi,t)$ and its mass points are the same as of
$d\rho_c(\xi),$ therefore \eqref{yudc} is hold for
$d\rho_c(\xi,t)$ as well. Thus, for each $t\in[0,T),$  for the
elements of \eqref{vl} corresponding to $d\rho_c(\xi,t),$ we have
that
\begin{equation*}
 0<\prod_{n=0}^{\infty} \frac{4}{c^2}a_n(t) < \infty,
\end{equation*}
and the claim follows.

Now consider the Volterra lattice \eqref{vl} in the infinite case.
As known, see \cite{ass,zh}, the elements of the Weyl matrix $\hat
\M(z,t)$ \eqref{wmp} satisfy the following relations for $Im\, z >
0 \,(Im\, z<0):$
\begin{eqnarray}
\label{mpm} m_{0,0}(z,t) =
\frac{m^-(z,t)}{1-a_{-1}(t)m^+(z,t)m^-(z,t)},
\; m_{1,1}(z,t) = \frac{m^+(z,t)}{1-a_{-1}(t)m^+(z,t)m^-(z,t)}\nonumber\\
m_{0,1}(z,t)=m_{1,0}(z,t)=
\frac{\sqrt{a_{-1}(t)}m^+(z,t)m^-(z,t)}{1-a_{-1}(t)m^+(z,t)m^-(z,t)},\qquad
\quad
\end{eqnarray}
where $m^+(z,t)$ and $m^-(z,t)$ are the Weyl functions \eqref{mvj}
of Jacobi operators ${\cal L}^+$ and ${\cal L}^-,$ generated by
the semi-infinite matrices
\begin{eqnarray*}
\quad L^+=L^+(t)=\begin{pmatrix} 0& \sqrt{a_0}& 0&0&0&\dots\\
           \sqrt{a_0}&0&\sqrt{a_1}&0&0&\dots\\
           0&\sqrt{a_1}&0&\sqrt{a_2}&0&\dots\\
            0&0&\ddots&\ddots&\ddots&\dots
            \end{pmatrix}\quad \text{and}\\ \qquad
            \\ L^-=L^-(t)=\begin{pmatrix} 0& \sqrt{a_{-2}}& 0&0&0&\dots\\
           \sqrt{a_{-2}}&0&\sqrt{a_{-3}}&0&0&\dots\\
           0&\sqrt{a_{-3}}&0&\sqrt{a_{-4}}&0&\dots\\
            0&0&\ddots&\ddots&\ddots&\dots
            \end{pmatrix}
\end{eqnarray*}
respectively. As follows from \eqref{mpm},
\begin{equation}
\label{detm}
\sqrt{a_{-1}(t)}(m_{0,0}(z,t)m_{1,1}(z,t)-m_{0,1}(z,t)^2)=m_{0,1}(z,t).
\end{equation}
Thus, if $m^+(z,t),$ $m^-(z,t)$ and $a_{-1}(t)$ are known, then
using \eqref{mpm} we can find the elements  of Weyl matrix \eqref
{wmp} and after that, applying the inverse Stieltjes transform,
reconstruct the spectral matrix defined in \eqref{spmatr}. Both
$m^+(z,t)$ and $m^+(z,t),$ are in turn the Stieltjes transforms
\eqref{mvj}-\eqref{mveven} of the spectral measures
$d\rho^+(\g,t)$ and $d\rho^-(\g,t)$ of the operators ${\cal L}^+$
and ${\cal L}^-$ respectively.

If the system \eqref{vl} with $n\in \mathbb{Z}$ has a finite
logarithmic Hamiltonian, and, additionally,
\begin{equation}
\label{trz} \sum_{n=-\infty}^{\infty} |a_n(t)-1|<\infty,
\end{equation}
then, applying Theorem \ref{th3}, we get the following result.
%%%%%%%%%%%%%%%%%%%%%%%%%%%%%%%%%%%%%%%%%%%%%
\begin{theorem}
\label{th4} To every solution of the infinite ($n\in \mathbb{Z}$)
system \eqref{vl} , satisfying \eqref{trz}, there correspond two
probability measures $d\rho^+(\g,t)$ and $d\rho^-(\g,t),$
satisfying for each $t\in [0,T)$ the same conditions as does the
measure $d\rho(\xi)$ of Theorem \ref{th3}, such that the triple
$\{d\rho^+(\g,t), d\rho^-(\g,t),a_{-1}(t)\}$ uniquely reconstructs
the spectral matrix \eqref{spmatr} of \eqref{vl}. Conversely, each
triple $\{d\rho^+(\g,t), d\rho^-(\g,t),f(t)\},$ where
$d\rho^\pm(\g,t)$ are the measures satisfying for each $t\in[0,T)$
the conditions of Theorem \ref{th3}, and $f(t)\in C^1[0,T),\,
f(t)>0$ corresponds to the spectral matrix $\hat \M$ of an
infinite system \eqref{vl} satisfying \eqref{trz}, with
$a_{-1}(t)=f(t)$.
\end{theorem}
%%%%%%%%%%%%%%%%%%%%%%%%%%%%%%%%%%%%%%%%%%%%
The last claim is partly based on the Theorem 1 of \cite{zh}
mentioned in Section 2. Namely, since both $d\rho^\pm (\g,t)$ have
infinitely many growth points, the corresponding Weyl functions
$m^\pm(z,t)$ \eqref{mvj} are not rational functions of $z$. Also,
if we define the matrix $m(z,t)=(m_{i,j}(z,t))_{i,j=0}^1$ as in
\eqref{mpm} with $a_{-1}(t)=f(t),$ we obviously  get the relation
\eqref{detm}, so all conditions of this theorem are fulfilled and
therefore $m(z,t)$ is the Weyl matrix \eqref{wmp} of a certain
operator $\hat {\cal L}$.

 Once the measures $d\rho^+(\g,t)$ and $d\rho^-(\g,t)$ are
known for the fixed $t$, we can calculate their moments, and then,
using \eqref{solisp}, find $(a_n(t))_{n=0}^\infty$ and
$(a_n(t))_{n=-2}^\infty$ respectively, which, alongside with
$a_{-1}(t),$  form the solution of \eqref{vl} in the infinite case
for such $t$. If we take $t=0,$ then $\{d\rho^+(\g,0),
d\rho^-(\g,0),a_{-1}(0)\}$ can be called the generating triple for
the lattice \eqref{vl}, \eqref{trz} but, unlike the semi-infinite
case, the measures $d\rho^+(\g,t)$ and $d\rho^-(\g,t)$ cannot be
found from \eqref{evrho} with $d\rho^+(\g,0)$ and $d\rho^-(\g,0)$
standing in place of $d\rho(\g,0)$. It follows from Theorem
\ref{th3}, that this situation corresponds to the case of
``broken" lattice \eqref{vl} with $a_{-1}(t)\equiv 0,$ when its
Weyl matrix \eqref{wmp} turns into
\begin{equation*}
 \begin{pmatrix} m^-(z,t) &
0\\
0 & m^+(z,t)
\end{pmatrix} \quad \text{with}\; m^\mp(z,t)=\frac
{d\rho^{\mp}(\g,t)}{z-\g}
\end{equation*}
 (see \eqref{mpm}) and such that both $\sum_{n=0}^\infty \ln
a_n(t)$ and $\sum _{n=2}^\infty \ln a_{-n}(t)$ are convergent, so
the logarithmic Hamiltonian $ H_0(t) = -\infty$ (also note that
since $a_n(t)\to 1$ as $n \to \pm \infty,$ we have $H_1(t)=\infty$
for the second Hamiltonian of \eqref{vl}). It remains an open and
interesting question in which situations the measures
$d\rho^+(\g,t)$ and $d\rho^-(\g,t)$ and their evolution in time
can be found explicitly.

Next consider the system
\begin{eqnarray}
\label{mvol}
\overset \cdot c_0 = -\frac{1}{c_1},\quad \overset
\cdot c_n =
\frac{1}{c_{n-1}}-\frac{1}{c_{n+1}},\quad n \ge 1. \\
c_n = c_n(t),\quad t \in [0,T), \quad 0<T\le \infty,\; \inf_{n\in
\mathbb{Z}_+} c_n(t)
> 0. \nonumber
\end{eqnarray}
In \cite{sur} a similar (up to the sign of the right-hand side)
system was called a modified Volterra lattice. In the coordinates
$(c_n)$ one can define the Poisson bracket $\{.,.\}_0$ with
non-vanishing elements
\begin{equation*}
\{c_{n+1},c_{n}\}_0 = -1,
\end{equation*}
\cite{sur} and write \eqref{mvol} in the Hamiltonian form
\begin{equation*}
\overset \cdot c_n = \{\hat H_0,c_{n}\}_0
\end{equation*}
with the Hamiltonian $\hat H_0 = \displaystyle \sum_{n=0}^\infty
\ln c_n$, so the question about the finiteness of the latter
remains actual. In studying this issue, we can use our previous
findings. Namely, assume that the system \eqref{mvol} has a
solution $(c_n(t))_{n=0}^\infty$ such that $\hat H_0$ is finite.
First, set $b_n = \displaystyle \frac{1}{c_n},$ then \eqref{mvol}
can be written in the equivalent form
\begin{equation}
\label{vbn}
\overset \cdot b_0 = b_0^2{b_1},\quad \overset \cdot
b_n = b_n^2(b_{n+1}-b_{n-1}),\quad n \ge 1.
\end{equation}
Next, by setting $a_n(t)=b_n(t)b_{n+1}(t),\; n \ge 0$ we get the
system \eqref{vl}-\eqref{seminf} with the initial data
\begin{equation}
\label{initabc}
(a_n(0)=b_{n}(0)b_{n+1}(0)=(c_{n}(0)c_{n+1}(0))^{-1})_{n=0}^{\infty}.
\end{equation}
As we see, in this case
\begin{equation}
\label{lnabc}
 \displaystyle \sum_{n=0}^\infty \ln b_n
= - \displaystyle \sum_{n=0}^\infty \ln c_n = \displaystyle
\frac{1}{2}(\sum_{n=0}^\infty \ln a_n + \ln b_0),
\end{equation}
 so $\hat H_0$ is
finite if so is the Hamiltonian $H_0$ of
\eqref{vl}-\eqref{seminf}. Also, if we assume that the solution of
\eqref{mvol} additionally satisfies the condition
\begin{equation}
\label{ctrace} \sum_{n=0}^{\infty} |c_n(t)-1|<\infty,
\end{equation}
then the same is true for $b_n(t)=1/c_n(t)$ and
$a_n(t)=b_n(t)b_{n+1}(t)$.
 Thus, summing up these observations
and using Theorem \ref{th3}, we come to the following conclusion
with respect to modified Volterra lattices.
%%%%%%%%%%%%%%%%%%%%%%%%%%%%%%%%%%%%%%
\begin{theorem}
\label{th5} To each solution of \eqref{mvol} satisfying
\eqref{ctrace} (and therefore having the finite Hamiltonian $\hat
H_0$) there corresponds the probability measure $d\rho(\g),$
satisfying the conditions of Theorem \ref{th3}, which coincides
with the spectral measure $d\rho(\g,0)$ of the operator ${\cal
L}(0)$ corresponding to the matrix $L(0)$ \eqref{lv} with the
entries defined from \eqref{initabc}.
\end{theorem}
%%%%%%%%%%%%%%%%%%%%%%%%%%%%%%%%%%%%%%%%%%%%%
Next, consider the Cauchy problem for \eqref{mvol} with the
initial data $(c_{n}(0))_{n=0}^\infty$. Note that this system
admits the Lax representation with non-Jacobi matrix $L$, see
\cite{osrjmp2} for more details, so the integration method
considered in Section 2 is not applicable here. Instead, using the
previous arguments, one can easily check that it can be reduced to
the one for the (classical) semi-infinite Volterra lattice.
Namely, we can first consider the Cauchy problem for
\eqref{vl}-\eqref{seminf} with the initial data
$(a_n(0)=(c_{n}(0)c_{n+1}(0))^{-1})_{n=0}^{\infty}$. Using the
algorithm of Theorem \ref{th1}, we find its solution
$(a_n(t)_{n=0}^{\infty}$ for $t\in [0,T)$. Next, we define the
functions $(b_n(t)_{n=0}^{\infty}$ as follows:
\begin{equation}
\label{ab}
\qquad b_0(t)=\displaystyle b_0(0)e^{\int_0^t
a_0(\tau)d\tau}, \quad b_n(t)=a_{n-1}(t)/b_{n-1}(t),\quad n \ge 1;
\end{equation}
which form the solution of \eqref{vbn} with the initial data
related to those of \eqref{vl}-\eqref{seminf} via \eqref{initabc};
and finally, by setting
\begin{equation}
\label{bc}
c_n(t)=1/b_n(t), \quad n \in \mathbb{Z}_+,
\end{equation}
we get the required solution of \eqref{mvol}. It follows from
Theorem \ref{th5} and \eqref{lnabc} that if the Hamiltonian $\hat
H_0$ of \eqref{mvol} is finite, then so is the one of the above
constructed system \eqref{vl}-\eqref{seminf}.

Conversely, we can start from a certain system
\eqref{vl}-\eqref{seminf} having $(a_n(0))_{n=0}^{\infty}$ as the
initial data and, consecutively using \eqref{ab} and \eqref{bc},
get the functions $(c_n(t))_{n=0}^{\infty}$ forming the solution
of \eqref{mvol} with the initial data satisfying \eqref{initabc}.
It is important to note that if the Hamiltonian $H_0(t)$ of
\eqref{vl}-\eqref{seminf} is finite and we define
$(b_n(t))_{n=0}^{\infty}$ as in \eqref{ab}, we cannot, generally
speaking, conclude that $H_b(t):=\sum_{n=0}^{\infty} \ln b_n(t),$
is finite as well (as follows from \eqref{bc}, the finiteness of
$H_b(t)$ and $\hat H_0(t)$ are equivalent to each other). Indeed,
if we set, e. g. $t=0$ and $a_n(0)=1$ for all $n,$ we obviously
get that $H_0(0)=0$ and, by setting $b_0(0)=1/2$ in \eqref{ab}, we
find that $b_n(0)=1/2$ for even $n$ and $b_n(0)=2$ when $n$ is
odd, so the series $H_b(0)$ is divergent. Since
$a_n(t)=b_n(t)b_{n+1}(t),$ we find that the partial sums
$S_N(H_b(t)):=\sum_{n=0}^N \ln b_n(t)$ of $H_b(t)$ satisfy the
following relations
\begin{equation*}
S_N(H_b(t)) = \frac{1}{2}(\ln b_0(t)+S_{N-1}(H_0(t))+\ln
b_{N+1}(t)),\quad N\ge 1,
\end{equation*}
where $S_{N-1}(H_0(t)):= \sum_{n=0}^{N-1} \ln a_n(t)$. Therefore,
$H_b(t)$ is finite if and only if \newline $\lim_{n\to \infty}
b_n(t)=1$ which is equivalent to
\begin{equation}
\label{limba} \lim_{n\to \infty} \frac{a_0(t)a_2(t)\dots
a_{2n}(t)}{b_0(t)a_1(t)\dots a_{2n+1}(t)}=1.
\end{equation}
Applying \eqref{solisp}, where $\Delta_n = |s_{i+j}(t)|_{i,j=0}^n$
and $s_{i+j}(t)$ are the moments of the spectral measure
$d\rho(\xi,t)$ corresponding to the system
\eqref{vl}-\eqref{seminf}, we find that \eqref{limba} can be
rewritten as follows
\begin{equation}
\label{limdel} \lim_{n\to \infty}
\frac{\Pi_1^n(t)}{b_0(t)\Pi_2^n(t)}=1,
\end{equation}
where
\begin{equation}
\label{pi12}
\Pi_1^n(t) = \frac{\Delta_1^2\Delta_3^2\dots
\Delta_{2n-1}^2 \Delta_{2n+1}}{\Delta_2^2 \Delta_4^2 \dots
\Delta_{2n}^2}, \quad \Pi_2^n(t)
=\frac{\Delta_{2(n+1)}}{\Pi_1^n(t)\Delta_{2n+1}}.
\end{equation}
Thus, substituting $a_0(t)=s_2(t)$ into \eqref{ab} and applying
Theorem \ref{th2}, we arrive at the following result.
\begin{theorem}
\label{th6}
 Each even probability measure $d\rho(\xi)$ with
essential support coinciding with the interval $[-2,2],$
satisfying there the the Szeg\H{o} condition and  having, possibly
infinite, point spectrum $\{\xi_j^{+}\},\,\{\xi_j^{-}\}, \;
\xi_j^{-} = - \xi_j^{+} $ satisfying the Peherstorfer-Yuditskii
condition \eqref{pu} generates the system \eqref{mvol} by
successively using the algorithm of Theorem \ref{th1} and
\eqref{ab}-\eqref{bc}. Such system has a finite Hamiltonian
$H_0=\sum_{n=0}^{\infty} \ln c_n(t)$ if the moments $s_k(t)$ of
the measure $d\rho(\xi,t)$ defined as in \eqref{rhohi} with
$d\tilde\rho(\xi)= d\rho(\xi),$ satisfy the following condition:

There exists $x>0$ such that
\begin{equation}
\label{xcond}
\lim_{n\to \infty} \frac{\Pi_1^n(t)}{e^{\int_0^t
s_2(\tau)d\tau}\Pi_2^n(t)}=x
\end{equation}
for each $t\in [0,T),$ where $\Pi_1^n(t)$ and $\Pi_2^n(t)$ are
defined in \eqref{pi12}
\end{theorem}
If \eqref{xcond} is fulfilled, then $x=b_0(0)$. Turning back to
Theorem \ref{th5}, we find that if \eqref{ctrace} is fulfilled,
then for the moments of the measure $d\rho(\g,t)$ defined in
\eqref{evrho} both \newline $\Pi_1(t)=\lim_{n\to\infty}
\Pi_1^n(t)$ and $\Pi_2(t)=\lim_{n\to\infty} \Pi_2^n(t)$ are
finite. Then it follows from \eqref{xcond} that
$x=\displaystyle\frac{\Pi_1(t)}{e^{\int_0^t
s_2(\tau)d\tau}\Pi_2(t)}$ is a first integral of \eqref{mvol}.
%Then, using  \eqref{evrho}-\eqref{momev} and \eqref{solisp}, one
%can find the solution $(a_n(t))_{n=0}^{\infty}$ of
%\eqref{vl}-\eqref{seminf} with the initial data The solutions
%$(b_n(t))_{n=0}^{\infty}$
% of \eqref{vbn} and,
%respectively, $(c_n(t))_{n=0}^{\infty}$ of \eqref{mvol} are then
%found as follows:

%Thus, if we have the probability measure satisfying the conditions
%of Theorem \ref{th3}, we can first generate the semi-infinite
%lattice \eqref{vl} with the finite logarithmic Hamiltonian, and
%then, using successively \eqref{ab} and \eqref{bc}, obtain the
%system \eqref{mvol} having the finite Hamiltonian $\hat H_0$.

%%%%%%%%%%%%%%%%%%%%%%%%%%%%%%%%%%%%%%%%%%%%%%%%%%%%%%%%%%%%%%%
%%%%%%%%%%%%%%%%%%%%%%%%%%%%%%%%%%%%%%%%%%%%%%%%%%%%%%%%%%%%%%%
\section{Conclusion. Open issues}
In view of the above, we can conclude that there exists a
correspondence between the even probability measures satisfying
the Szeg\H{o} and Peherstorfer-Yuditskii conditions and the
Volterra lattices with a finite logarithmic Hamiltonian. The
Szeg\H{o} condition can be interpreted as finiteness of measure's
entropy, see \eqref{segtrig}, while the Peherstorfer-Yuditskii
condition allows its mass point part to be infinite. It is of
interest to give a proper physical interpretation of our findings,
especially taking into account that \eqref{vft} provides a link
between the Volterra lattice and the conformal field  theory
\cite{ft1,fv}.

Turning back to Theorem 3, we have to conclude that finding a
description of the subclass of probability measures, satisfying
the conditions of this theorem, which one-to-one corresponds to
the class of semi-infinite Volterra lattices \eqref{vl} such that
\newline $\sum_{n=0}^{\infty} |a_n(t)-1|<\infty,$ remains an
unresolved issue. We believe that the key step here will be to
establish a criterion of resolution of the inverse spectral
problem for Jacobi operators of the trace class. Note that for the
Hilbert-Schmidt class, such criterion have been obtained in the
Theorem 1 of the above-mentioned paper \cite{kils} by Killip and
Simon (see also Theorem 13.8.6 of \cite{simop}).The latter result
can be regarded as a contribution to the theory of incomplete
inverse problems, i. e. the ones where part of the information
about the object to reconstruct is known in advance, see
\cite{yurm,osrjmp} for more details, so such activity is also
important for this theory. In the context of our Theorem 1 and
\eqref{marc}, it remains an actual task to study in which cases
the point spectrum of a Jacobi operator is finite.

As mentioned in the Introduction, the Bogoyavlensky lattices, e.
g. the ones satisfying the equations
\begin{equation}
\label{bogl} \overset \cdot a_n = a_n\left(\sum_{k=1}^q
a_{n+k}-\sum_{k=1}^q a_{n-k}\right)
\end{equation}
for some fixed $q\ge 1,$ can be regarded as Voltera-type systems.
In the semi-infinite case, their integration methods by means of
the inverse spectral problems for band operators were considered
in \cite{akv,osajm}. The systems \eqref{bogl} inherit many
properties of the Volterra lattice (e. g. they can be regarded as
discrete versions of the Korteweg - de Vries equation
\cite{bogo}), so the applicability of our findings to their case
is an interesting and difficult task. In this connection, the
paper \cite{aks} of Aptekarev, Kaliagin and Saff, where, in
particular, the spectral data (elements of the Weyl function) of
the sparse band operators of the trace class, i. e. such  that
$\sum_{n=0}^{\infty} |a_n-1|<\infty,$ have been studied, is of
certain interest. Also, a natural step in this direction is to
find an appropriate Poisson algebra, by using which one can get a
representation for \eqref{bogl} similar to \eqref{vft}. To our
knowledge, see \cite{sur} and references thereafter, this has not
yet done.

All these issues (including the open ones mentioned in the
previous Section) can be the subject of further research.

This work is done at SRISA according to the project FNEF-2024-0001
(Reg. No 1023032100070-3-1.2.1).

%%%%%%%%%%%%%%%%%%%%%%%%%%%%%%%%%%%%%%%%%%%%%%%%%%%%%%%%%%%%%
%%%%%%%%%%%%%%%%%%%%%%%%%%%%%%%%%%%%%%%%%%%%%%%%%%%%%%%%%%%%%%%%%

\end{document}